\documentclass[12pt]{article}
\usepackage{amsmath,amsthm,amssymb}

\renewcommand{\theequation}{\mbox{\arabic{section}.\arabic{equation}}}

\newtheorem{theorem}{Theorem}[section] 
\newtheorem{lemma}{Lemma}[section]

\newtheorem{ass}{Assumption}[section]
\newtheorem{observe}{Observation}[section]
\newtheorem{proposition}{Proposition}[section]

\numberwithin{equation}{section}

\newcommand{\D}{{\rm d}}
\newcommand{\dx}{\, \D x}
\newcommand{\dy}{\, \D y}

\newcommand{\dis}{\displaystyle}

\newcommand{\msp}{\;\;}
\newcommand{\fsp}{\quad\;}
\newcommand{\psp}{\,}

\newcommand{\diff}{\mathcal{D}}

\newcommand{\rz}{\mathbb{R}}
\newcommand{\nz}{\mathbb{N}}

\newcommand{\iom}{\int_{\Omega}}

\newcommand{\klauf}{\left(\begin{array}}
\newcommand{\klzu}{\end{array}\right)}

\title{Existence results for generalized EED denoising problems}
\author{Michael Bildhauer \& Martin Fuchs}
\date{}

\newcommand{\br}{\mathcal{B}_R}

\usepackage{hyperref}

\hypersetup{
   pdfnewwindow=true,
   colorlinks=false,
   linkcolor=black,
   citecolor=green,
   filecolor=magenta,
   urlcolor=cyan,
}

\begin{document}
\parindent0em
\maketitle

\newcommand{\hypref}[2]{\hyperref[#2]{#1 \ref*{#2}}}
\newcommand{\hypreff}[1]{\hyperref[#1]{(\ref*{#1})}}

\begin{abstract}
The joint work of the authors with Marcelo C\'ardenas and Joachim Weickert \cite{BCFW:2019_1} 
on edge-enhancing diffusion inpainting problems leads to the analysis of 
related denoising problems.
Here, a surprisingly broad class of diffusion tensors is
admissible to obtain the existence of solutions to EED denoising problems.\footnote{The
authors thank Joachim Weickert for a series of stimulating discussions which
also inspire a forthcoming numerical analysis of the problem under consideration.}
\end{abstract}

\section{Introduction and main result}\label{intro}

One possibility of restoring missing image data is given by the well known EED-Ansatz 
(edge-enhancing diffusion).
This kind of inpainting problem is discussed,
for instance, in \cite{BCFW:2019_1}, where the reader will find the necessary background material
including a list of references.\\

Roughly speaking, given an ''image-data set'' $K\subset \Omega$ and a function $u\in L^1(G)$, $G:=\Omega- K$,
the mollification $u_\sigma$ w.r.t.~$G$ and with Gaussian kernel is considered as argument of the
diffusion tensor $\overline{D}$.\\ 

In principle, $\overline{D}$ looks like (in fact its representation in the tangential 
and in the normal direction of the edges, respectively)
\begin{equation}\label{intro 1}
\overline{D}\big(\nabla (u_\sigma)\big) := 
\klauf{cc} 1&0\\[1ex]0&\big(1+|\nabla (u_\sigma)|^2\big)^{-1/2}\klzu \psp .
\end{equation}

The problem then reads as: Find a solution $u$ of the boundary value problem
\begin{eqnarray}
\label{intro 2}
{\rm div}\psp \Big(\overline{D}\big(\nabla (u_\sigma)\big)\nabla u\Big)\psp =\psp  0&\mbox{on}&G\psp ,\\[1ex]
\label{intro 3}
u\psp = \psp f&\mbox{on}& \partial K\psp ,\\[1ex]
\label{intro 4}
\overline{D}\big(\nabla (u_\sigma)\big)\nabla u \cdot \mathcal{N} \psp = \psp 0&\mbox{on}& \partial \Omega\psp ,
\end{eqnarray}
where $\mathcal{N}$ denotes the outward unit normal to $\Omega$ and $f$: $K \to \mathbb{R}$
is a given function representing the image data.\\

In \cite{BCFW:2019_1}, the existence of a solution to \hypreff{intro 2} -- \hypreff{intro 4} is proved  
using a Leray-Schauder fixed point argument. Moreover, some additional analytical results on the 
fixed point set and a priori estimates for particular iterations
are established.\\

Note that the arguments outlined in \cite{BCFW:2019_1} are strongly adapted to
the limit case $\mu=1$ for the second eigenvalue
\[
(1+|\nabla (u_\sigma)|^2)^{-\mu/2}
\] 
of the diffusion tensor in \hypreff{intro 1}. We also note that this limit case corresponds to linear
growth as lower bound for $\overline{D}(Z)Z$.\\

The purpose of our considerations below is to replace the above described EED inpainting problem
by a wide class of EED-related denoising models
which surprisingly allows us to include variants with very weak assumptions on the
generalized diffusion tensor, for example, any value of $\mu$ is admissible.\\

We first observe that the new problem differs in the structure of the underlying equations in the sense
that \hypreff{intro 2} should be valid on 
the whole domain $\Omega$ with a non-vanishing suitable r.h.s. The condition
\hypreff{intro 3} is omitted and \hypreff{intro 4} is carried over.\\

To be precise, let us first have a closer look at denoising procedures from the abstract point of view.
The interested reader is referred to the monograph \cite{we}.

\begin{enumerate}
\item \emph{Data term for denoising problems.}

Let $t \geq 1$ be some fixed number in the following and suppose that we are given data (observed image)
$f$ of class $L^t(K)$ where, depending on the problem under consideration,  $K = \Omega$ or 
$K$ is a suitable subset of $\Omega$.\\

As a measure for data-fitting we may consider any strictly convex functional
\[
E[K,\cdot,f]: \psp L^t(K) \to \rz^+_0\psp ,\fsp v \mapsto \int_K e(v,f) \dx\psp ,
\]
where the density $e(\cdot,f)$ is a smooth and strictly convex function.\\
 
Of course the most prominent example is given with the choices $K=\Omega$, $t=2$ and ($\lambda >0$ fixed)
\[
E[v]= E[\Omega,v,f] = \lambda \iom |v-f|^2 \dx \psp .
\]

Having our main applications discussed below in mind, we suppose throughout this paper
that the data term is given by ($K$, $f$ fixed as above)
\[
E[K,v,f] = \lambda \int_K |v-f|^t \dx \psp ,\fsp t > 1\psp .
\]

\item \emph{Regularization in variational form.} (compare \cite{BF:2012_1} and 
the references quoted therein)

Given a particular convex energy density $R$: $\rz^2\to \rz$, the functional
\[
w \mapsto \iom R(\nabla w) \dx
\]
is considered as regularizing part defined on an appropriate function space
and from the variational point of view we obtain
\[
\iom \nabla R(\nabla u)\cdot \nabla \varphi \dx \fsp\mbox{for all admissible}\msp
\varphi :\psp \Omega \to \rz
\]
as the leading contribution to the Euler equation. In the strong sense this is written as
\[
-{\rm div} \big[ \nabla R(\nabla u)\big]\psp .
\]

\item \emph{EED regularization.} (compare \cite{BCFW:2019_1})

If $u_\sigma$ denotes the mollification of $u$ via a kernel $k_\sigma$, $\sigma > 0$ fixed,
then we may replace $\nabla R(\nabla u)$ by the diffusion term
\begin{equation}\label{intro 5}
-{\rm div}\psp \big[\overline{D}\big(\nabla (u_\sigma)\big) \nabla u\big] 
\end{equation}
with suitable Eigenvalues of $\overline{D}$ depending on the edge direction.\\
\end{enumerate}

In our note we generalize \hypreff{intro 5} by the way admitting
a large amount of flexibility in choosing appropriate regularizations.\\

At this point we may formulate the general assumption of our considerations:

\begin{ass}\label{intro ass 1}
Suppose that $\Omega \subset  \rz^2$ is a bounded Lipschitz domain, $K\subset \Omega$ is a set
of positive measure, $f$ is of class $L^t(K)$ for some $t>1$ and $\lambda >0$ is some arbitrary fixed
parameter. \\

For any $R>0$ and $t$ fixed as above we define the ``data ball''
\[
\br^K := \{w \in W^{1,2}(\Omega): \psp \|w\|_{L^t(K)} \leq R\}
\]
being a closed and convex subset of the space $W^{1,2}(\Omega)$.
Let
\begin{eqnarray*}
\rz^{2\times 2}_{sym,+} &:=& \Big\{ (F_{\alpha \beta})_{1\leq \alpha,\beta \leq 2}: \psp
F_{\alpha\beta} =F_{\beta \alpha}, \\
&&\fsp F Q\cdot Q > 0 \msp \mbox{for all}\msp Q\in \rz^2,\psp Q\not=0\Big\}
\end{eqnarray*}
and suppose that the generalized diffusion tensor $\diff$ satisfies
\begin{equation}\label{intro 6}
\diff:\psp W^{1,2}(\Omega) \to C^0\big(\overline{\Omega},\rz^{2\times 2}_{sym,+}\big)\msp\mbox{is a continuous mapping}\psp .
\end{equation}
As a final hypothesis we take:
\begin{equation}\label{intro 7}
\mbox{For all $R>0$ the set $\diff(\br^K)$ is precompact in $C^0\big(\overline{\Omega},\rz^{2\times 2}_{sym,+}\big)$.}
\end{equation}
Referring 
to the theorem of Arzela and Ascoli, we observe that \hypreff{intro 7} is a consequence of the 
stronger condition
\addtocounter{equation}{-1}
\renewcommand{\theequation}{\mbox{\arabic{section}.\arabic{equation}$^*$}}
\begin{eqnarray}
\nonumber &\mbox{For  some $\alpha \in (0,1]$ and for all $R>0$}&\\[1ex]
\label{intro 3ast}
&\mbox{the set 
$\diff(\br^K)$ is a bounded subset of $\dis C^{0,\alpha}\big(\overline{\Omega},\rz^{2\times 2}_{sym,+}\big)$.}&
\end{eqnarray}
\renewcommand{\theequation}{\mbox{\arabic{section}.\arabic{equation}}}
\end{ass}

Let us formulate our main result.
\begin{theorem}\label{intro theo 1}
If \hypref{Assumption}{intro ass 1} holds true, then there exists a function $u\in W^{1,2}(\Omega)$ such that 
\begin{eqnarray}\label{intro 8}
{\rm div}\psp \Big[ \diff(u)\nabla u\Big] &=& \lambda \psp \chi_K \psp t |u-f|^{t-2} (u-f)\fsp\mbox{in $\Omega$}\psp ,\\[2ex]
\label{intro 9}
\diff(u) \nabla u \cdot \mathcal{N} &=& 0 \fsp\mbox{on $\partial \Omega$}\psp .
\end{eqnarray} 

Here $\chi_K$ denotes the characteristic function of the set $K$ and 
we use the symbol $\mathcal{N}$ for the outward unit normal to $\partial \Omega$.\\
\end{theorem}

In order to prove \hypref{Theorem}{intro theo 1}, we define a related operator
$T$ in \hypref{Section}{op} and establish some fundamental properties in \hypref{Section}{T}.\\

Using well known fixed point arguments applied to the operator $T$
we obtain \hypref{Theorem}{op theo 1}
by the way proving \hypref{Theorem}{intro theo 1} as a corollary.\\

Having these preparations and notation in mind, we now like to present a short
section with a list of examples.\\

In case of the examples i) - iii) below, it is quite easy to check \hypref{Assumption}{intro ass 1}. 
Concerning the most interesting examples iv) and v), 
we postpone the proof to \hypref{Section}{inpaint}.

\section{Some examples for generalized diffusion tensors}\label{diff}

\begin{enumerate}
\item We first discuss tensors $\diff$ defined via a smoothing procedure by the way extending
example \hypreff{intro 1} from the introduction.

We take $K=\Omega$ and let (with Gaussian kernel $k_\sigma$, $\sigma > 0$ fixed)
\begin{equation}\label{diff 1}
w_\sigma (x) := \iom k_\sigma(y-x) w(y)\dy \fsp\mbox{for any}\msp w\in L^1(\Omega)\psp .
\end{equation}
We then let for $k\in \nz_0$
\[
\diff_{\sigma}^{(k)}(u):= \Big(d_{\alpha \beta}(\cdot,u_\sigma, \nabla (u_\sigma),\dots , \nabla^k(u_\sigma)
\Big)_{1\leq \alpha, \beta \leq 2}
\]
with given continuous coefficients
\[
d_{\alpha \beta}:\psp \overline{\Omega}\times \rz \times \rz^2\times\dots \times \rz^{2k}\to \rz \psp ,\fsp
d_{\alpha \beta} =d_{\beta \alpha} \psp ,
\]
such that 
\[
\sum_{\alpha, \beta =1}^{2} d_{\alpha \beta}(\dots)q_\alpha q_\beta > 0
\]
holds for all $q\in \rz^2-\{0\}$ and any argument of $d_{\alpha \beta}$.

A simple explicit example in the spirit of \hypreff{intro 1} is given by
\[
\diff_1(u)  := 
\klauf{cc} 1&0\\[1ex]0&\big(1+|\nabla (u_\sigma)|^2\big)^{-\mu/2}\klzu \fsp\mbox{for some}\msp \mu \in \rz\psp ,
\]
which satisfies \hypref{Assumption}{intro ass 1} as it is shown in \cite{BCFW:2019_1}.

\item Motivated by the $p$-Laplacian, we would like to incorporate a 
''weighted'' diffusion like $\diff_1(u) |\nabla u|^p$, $p>1$, in our considerations.

Although this is not obvious, we find a good approximation by introducing a second smoothing
parameter $\delta >0$ via \hypreff{diff 1}:
\[
\diff_2(u)  := 
\klauf{cc} 1&0\\[1ex]0&\big(1+|\nabla (u_\sigma)|^2\big)^{-\mu/2}\klzu |\nabla (u_\delta)|^p 
\]
(or more general $\diff_{\sigma}^{(1)}(u)|\nabla (u_\delta)|^p$).

Notice that both in example i) and example ii) we use the fact that the data term is defined w.r.t.~the
whole domain $\Omega$. 
Moreover, note that \hypref{Assumption}{intro ass 1} holds for any growth rate $t$.

\item A completely different type of example is given in the spirit of ''Galerkin-type'' methods. 

Suppose that $t=2$, $K=\Omega$ and that we have fixed some functions 
$v_k \in C^{\infty}_0(\Omega)$, $k=1,$ \dots, $N$.\\

We then consider diffusion tensors of type
\[
\diff_3(u) = \tilde{\diff}_3 \Bigg(\sum_{k=1}^{N} (u,v_k)_{L^2(\Omega)} \psp v_k\Bigg)
\]
with appropriate choice of $\tilde{\diff}_3$ s.t.~\hypref{Assumption}{intro ass 1} holds.

\item In the next example we use some kind of preconditioning with a standard denoising in order
to map the data ball in a set of smooth functions in the sense of \hypref{Assumption}{intro ass 1}.

We consider the case $K=\Omega$ and we suppose that $t>2$.
Given $w \in \br^\Omega$ we start with a Whittaker-Tikhonov regularization 
of the function $w$, i.e.~we first denoise $w$ by solving
\[
\iom |\nabla v|^2 \dx + \iom |v-w|^2\dx \to \min \fsp\mbox{in}\msp W^{1,2}_0(\Omega) \psp .
\]
In \hypref{Theorem}{diff theo 1} below we will show that
the solution $\hat{v} = \hat{v}(w)$ is bounded in $C^{1,\alpha}(\overline{\Omega})$ (uniformly w.r.t.~$w$),
hence we may consider the tensor
\[
\diff_4 (u) = \tilde{\diff}_4 \big(\hat{v}(u)\big) \psp ,
\]
where for example
\[
\tilde{\diff}_4(w):= \Big(d_{\alpha \beta}(\cdot,w,\nabla w)\big)_{1\leq \alpha,\beta \leq 2}
\]
for a suitable matrix $(d_{\alpha\beta}(x,y,p))_{1\leq a, \beta \leq 2}$.

Keeping this example in mind, we finally include the case $K\subset \Omega$ in the next example
and formulate \hypref{Therorem}{diff theo 1} to verify our assumption.

\item Suppose that we are given continuous functions $d_{\alpha \beta}$, 
$d_{\alpha \beta} = d_{\beta \alpha}$, $1\leq \alpha$, $\beta \leq 2$,
\[
\big(d_{\alpha \beta}(x,y,z)\big)_{1\leq \alpha,\beta\leq 2}:\psp \overline{\Omega} \times \rz \times \rz^2 \to \rz \psp .
\]

For $t > 1$ and a measurable set $K \subset \Omega$ of positive measure  
we fix $1 < s$ and define the inpainting operator
\[
I:\psp W^{1,2}(\Omega)\to W^{1,2}_0(\Omega)\psp , \fsp w \mapsto u=I(w) \psp ,
\]
where $u$ is the unique solution of the minimization problem ($\lambda > 0$ fixed)
\begin{equation}\label{diff 2}
\iom |\nabla u|^2 \dx + \lambda \int_K |u-w|^s\dx \to \min \fsp \mbox{in}\msp W^{1,2}_0(\Omega)\psp .
\end{equation}
\end{enumerate}

We have \hypref{Theorem}{diff theo 1}, which will be proved in \hypref{Section}{inpaint}.

\begin{theorem}\label{diff theo 1}
Suppose we are given the inpainting operator $I$ of the last example. With the choice
\begin{equation}\label{diff 3}
1 < s < 1+ \frac{t}{2}
\end{equation}
we let
\[
\diff(w):= \Big(d_{\alpha \beta}\big(\cdot, I(w), \nabla I(w)\big)\Big)_{1\leq \alpha ,\beta \leq 2} \psp .
\]
Then \hypref{Assumption}{intro ass 1} holds true, provided $\Omega$ is a $C^{1,1}$ domain.
\end{theorem}

\section{The operator $T$}\label{op}

We consider the functional
\begin{eqnarray*}
J[\cdot, \cdot,\cdot]:&&W^{1,2}(\Omega)\times W^{1,2}(\Omega)\times L^t(K) \to \rz \psp ,\\[2ex]
J[w,v,f]&:=&\frac{1}{2} \iom \diff(w) \nabla v \cdot \nabla v\dx + \lambda  \int_K |v-f|^t \dx \psp ,
\end{eqnarray*}
and the minimization problem
\begin{equation}\label{op 1}
J[w,\cdot,f]\to \min\fsp \mbox{in}\msp W^{1,2}(\Omega) \psp .
\end{equation}

We define ($f$ fixed) the operator $T$: $W^{1,2}(\Omega) \to W^{1,2}(\Omega)$,
\[
w \mapsto u:= T(w) \psp ,
\]
where $T(w)$ denotes the unique solution of problem \hypreff{op 1}.\\

The Euler equation for the problem under consideration reads as
\begin{equation}\label{op 2}
\iom \diff(w) \nabla u \cdot \nabla \varphi\dx + \lambda \int_K t |u-f|^{t-2} (u-f) \varphi \dx  = 0 
\end{equation}
for all $\varphi\in W^{1,2}(\Omega)$.
Note that, if a fixed point $u$ of $T$ is found, then we have a weak solution of \hypreff{intro 8}.\\

Another essential tool for our considerations is\

\begin{observe}\label{op observe 1}
Suppose that $R$ is sufficiently large. Then we have
\begin{equation}\label{op 3}
T\big(W^{1,2}(\Omega) \big) \subset \br^K \psp ,\fsp\mbox{in particular}\fsp
T\big( \br^K \big) \subset \br^K \psp . 
\end{equation}
\end{observe}

In fact, given $R$ sufficiently large, the claim follows from 
\[
E[K,T(w),f] \leq J[w,T(w),f] \leq J[w,T(0),f] = E[K,0,f]\leq R^{t}\psp .
\]

At this point we formulate our main result:

\begin{theorem}\label{op theo 1}
If $R>0$ is sufficiently large, then the mapping 
$T$: $W^{1,2}(\Omega) \to W^{1,2}(\Omega)$ has at least one fixed point 
in $\br^K$.
\end{theorem}

\emph{Proof.} With \hypreff{op 3} of Observation \ref{op observe 1}, \hypref{Lemma}{T lem 1}
and \hypref{Lemma}{T lem 2} of the next section, 
we have verified all hypotheses needed for Corollary 11.2 of \cite{GT}.\hspace*{\fill}$\Box$\\

\section{Properties of $T$}\label{T}

\begin{lemma}\label{T  lem 1}
$T$ is a continuous operator $W^{1,2}(\Omega) \to W^{1,2}(\Omega)$.
\end{lemma}

\emph{Proof.} We conder a sequence $\{w_n\}$ and a function $w$ from $W^{1,2}(\Omega)$
such that as $n\to \infty$

\begin{equation}\label{T 1}
\|w_n-w\|_{W^{1,2}(\Omega)} \to 0
\end{equation}
and letting $u_n:= T(w_n)$, $u=T(w)$, we claim
\begin{equation}\label{T 2}
\|u_n-u\|_{W^{1,2}(\Omega)} \to 0 \, ,
\end{equation}
again as $n\to \infty$.\\

To this purpose we first establish the uniform bound
\begin{equation}\label{T 3}
\sup_{n\in \nz} \|u_n\|_{W^{1,2}(\Omega)} < \infty \psp .
\end{equation}

In fact, $\diff(w_n)(x)$ is positive definite for any $x\in \overline{\Omega}$, we have \hypreff{intro 6} and the   
sequence $\{w_n\}$ is bounded in $\br^K$ according to \hypreff{T 1}. 
Thus, the minimality of $u_n$ gives
\begin{eqnarray*}
\iom |\nabla u_n|^2\dx & \leq & c \iom \diff(w_n) \nabla u_n \cdot \nabla u_n\dx\\[2ex]  
&\leq & c J[w_n,u_n,f] \leq c J[w_n,0,f] \leq c\psp .
\end{eqnarray*}

Now, using \hypreff{T 3}, we may consider a subsequence $(u_{n_k})$ s.t.
\begin{equation}\label{T 4}
u_{n_k} \rightharpoondown : \tilde{u}\fsp\mbox{in}\msp W^{1,2}(\Omega)\msp\mbox{as}\msp k\to \infty \psp .
\end{equation}

Let us have a closer look at \hypreff{op 2} w.r.t.~$w_n$ and $u_n$: we have
\begin{equation}\label{T 5}
\iom D(w_n)  \nabla u_n  \cdot \nabla \varphi \dx = \lambda t \int_K |u_n-f|^{t-2} (u_n-f) \varphi \dx \psp .
\end{equation}
We benefit from \hypreff{intro 7} and, using \hypreff{T 4}, we pass to the limit in \hypreff{T 5} to obtain
\begin{equation}\label{T 6}
\iom D(w)  \nabla \tilde{u} \cdot \nabla \varphi dx =  \lambda t \int_K |\tilde{u}-f|^{t-2} (\tilde{u}-f) \varphi \dx\psp ,
\end{equation}
which by the uniqueness of solutions to \hypreff{T 6} implies $\tilde{u} = T(w) =u$.\\

This holds for any convergent subsequence, hence \hypreff{T 4} is true for the whole sequence
with $u=\tilde{u}$.\\

We now have to improve \hypreff{T 4} in the sense
\begin{equation}\label{T 7}
\nabla u_{n} \to \nabla u\fsp\mbox{in}\msp L^{2}(\Omega,\rz^2)\msp\mbox{as}\msp k\to \infty \psp .
\end{equation}

For proving \hypreff{T 7} we make use of \hypreff{T 5}, \hypreff{T 6} and observe ($\varphi := u_n - u$)
\begin{eqnarray*}
\lefteqn{\iom |\nabla u_n - \nabla u|^2 \dx
\leq  c \iom \diff(w) (\nabla u_n - \nabla u)
\cdot (\nabla u_n - \nabla u) \dx}\\[2ex]
&=& c\iom \diff(w) \nabla u_n \cdot (\nabla u_n - \nabla u)\dx
- c \iom \diff(w) \nabla u \cdot (\nabla u_n - \nabla u)\dx\\[2ex]
&=&c  \iom \big[ \diff(w) -  \diff(w_n)\big] \nabla u_n \cdot (\nabla u_n - \nabla u)\dx\\[2ex]
&& +c \iom \diff(w_n) \nabla u_n \cdot (\nabla u_n - \nabla u)\dx 
+ c\psp  \lambda t \int_K |u-f|^{t-2} (u-f) \varphi \dx\\[2ex]
&=&\iom \big[\diff(w) - \diff(w_n)\big] \nabla u_n \cdot (\nabla u_n - \nabla u)\dx\\[2ex]
&& + c \lambda t  \Bigg[\int_K |u-f|^{t-2} (u-f) \varphi \dx - \int_K |u_n-f|^{t-2} (u_n-f) \varphi \dx \Bigg]\psp .
\end{eqnarray*}
Again, by \hypreff{intro 7}, the smallness of the first term on the r.h.s.~is evident,
the second converges to zero on account of the weak $W^{1,2}$-convergence.
\hspace*{\fill}$\Box$\\

\begin{lemma}\label{T lem 2}
The set $T\big(\br^K\big)$ is precompact.
\end{lemma}

\emph{Proof.} Consider a sequence $(u_n) \subset T(\br^K)$, $u_n= T(w_n)$ for some $w_n \in \br^K$.\\

We claim that we can extract a subsequence, which is strongly converging in $W^{1,2}(\Omega)$.\\

As above, we have the uniform bound
\begin{equation}\label{T 8}
\sup_{n}\|u_n\|_{W^{1,2}(\Omega)} < \infty 
\end{equation}
and passing to a subsequence we may suppose
\begin{equation}\label{T 9}
u_{n_k} \rightharpoondown: \tilde{u}\fsp\mbox{in}\msp W^{1,2}(\Omega)\psp ,
\fsp u_{n_k} \to \tilde{u} \fsp\mbox{in}\msp L^{s}(\Omega)\msp\mbox{for any $1< s <\infty$}\psp .
\end{equation}

Now let
\[
A_{n_k}:\psp \overline{\Omega} \to \rz^{2\times 2}_{sym,+}\psp ,\fsp  A_{n_k}(x) := \diff(w_{n_k})(x) \psp .
\]
We have $w_{n_k} \in \br^K$, hence by \hypreff{intro 7} we may suppose the uniform convergence
to a function $A$:
\begin{equation}\label{T 10}
A_{n_k} \rightrightarrows A:\psp \overline{\Omega}\to \rz^{2\times 2}_{sym,+}\psp .
\end{equation}
Observe the estimate
\begin{eqnarray}\label{T 11}
\lefteqn{\iom A(x) (\nabla u_{n_k} - \nabla \tilde{u}) \cdot (\nabla u_{n_k} - \nabla \tilde{u})\dx}\nonumber\\[2ex]
&=& \iom A_{n_k}(x) (\nabla u_{n_k} - \nabla \tilde{u}) \cdot (\nabla u_{n_k} - \nabla \tilde{u})\dx\nonumber\\[2ex]
&& + \iom \big[A(x)-A_{n_k}(x)\big] (\nabla u_{n_k} - \nabla \tilde{u})
\cdot (\nabla u_{n_k} - \nabla \tilde{u})\dx \nonumber\\[2ex]
&=:& \alpha_k + \beta_k \psp ,
\end{eqnarray}
where \hypreff{T 8} and \hypreff{T 10} immediately give
$\lim_{k\to \infty} \beta_k = 0$.
Moreover,
\begin{eqnarray}\label {T 12}
\alpha_k &=& \iom A_{n_k}(x) \nabla u_{n_k} \cdot (\nabla u_{n_k} - \nabla\tilde{u}) \dx
- \iom A_{n_k}(x) \nabla \tilde{u} \cdot (\nabla u_{n_k} - \nabla\tilde{u}) \dx\nonumber\\[2ex]
&=& - \iom A_{n_k}(x) \nabla \tilde{u} \cdot (\nabla u_{n_k} - \nabla\tilde{u}) \dx \nonumber\\[2ex]
&&- \lambda t \int_K |u_{n_k} -f|^{t-2} (u_{n_k} -f) (u_{n_k} - \tilde{u}) \dx
\end{eqnarray}
since we have the Euler equation for $u_{n_k} = T(w_{n_k})$.\\

Again \hypreff{T 10} and weak convergence yield the convergence of the first term in \hypreff{T 12},
the second one is handled with the strong convergence stated in \hypreff{T 9}, thus
\[
\lim_{k\to \infty} \alpha_k = 0 
\]
and \hypreff{T 11} shows 
\[
\lim_{k\to \infty} \iom A(x) (\nabla u_{n_k} - \nabla \tilde{u}) \cdot (\nabla u_{n_k} - \nabla \tilde{u})\dx = 0 \psp .
\]
A final application of \hypreff{intro 6} and \hypreff{intro 7} leads to
\[
\nabla u_{n_k} \to \nabla \tilde{u} \fsp\mbox{in}\msp L^2(\Omega,\rz^2) \psp .
\]
and we have found a sequence strongly converging in $W^{1,2}(\Omega)$ which completes the proof
of \hypref{Lemma}{T lem 2}.
\hspace*{\fill}$\Box$\\

\section{Proof of \hypref{Theorem}{diff theo 1}}\label{inpaint}
In this last section we are going to prove \hypref{Theorem}{diff theo 1} and we always
refer to the operator $I$ as defined in Example v) of \hypref{Section}{diff}.
We also use the notation introduced in \hypref{Assumption}{intro ass 1}.\\

We start by establishing the compactness of $I$.

\begin{proposition}\label{inpaint prop 1}
Suppose that $\partial \Omega$ is of class $C^{1,1}$ and that (recalling \hypreff{diff 3})
\begin{equation}\label{inpaint 1}
1 < s < \frac{t}{2} +1 \psp ,\fsp \mbox{i.e.}\msp p:= \frac{t}{s-1} > 2 \psp .
\end{equation}
Then the set $I(\br^K)$ is bounded in $W^{2,p}(\Omega)$, hence bounded in 
$C^{1,\alpha}(\overline{\Omega})$ choosing $\alpha = (p-2)/2$.
\end{proposition}
\emph{Proof.} Fix $w\in \br^K$ and let $u =I(w)$ as defined in \hypref{Section}{diff}, Example v), i.e.~$u$
is the solution of the minimization problem \hypreff{diff 2} which gives using \hypreff{inpaint 1}
\[
\iom |\nabla u|^2 \dx \leq c \int_K |w|^s \dx \leq c \int_K |w|^t \dx \psp .
\]
As an immediate consequence of Sobolev's inequality we obtain the uniform bound
\begin{equation}\label{inpaint 2}
\|u\|_{L^{q}(\Omega)} \leq c(q,R) \fsp\mbox{for all}\msp q < \infty \psp .
\end{equation}
Moreover, $u$ satisfies for all $\varphi\in C^{\infty}_{0}(\Omega)$
\[
\iom \nabla u \cdot \nabla \varphi \dx + \int_K  s |u-w|^{s-2}(u-w) \varphi \dx = 0 \psp ,
\]
i.e.~we have in the weak sense
\begin{eqnarray}\label{inpaint 3}
- \Delta u&=& g \fsp\mbox{on}\msp\Omega\psp ,\fsp u=0\msp\mbox{on $\partial \Omega$}\psp ,\\[2ex]
\label{inpaint 4} g&:=& - s |u-w|^{s-2} (u-w) \chi_K  \psp .
\end{eqnarray}
Note that by \hypreff{inpaint 2} and on account of $w\in \br^K$ we have uniformly
\begin{equation}\label{inpaint 5}
\|g\|_{L^{t/(s-1)}(\Omega)} = \|g\|_{L^p(\Omega)}  \leq c(R)\psp .
\end{equation}

Referring to Theorem 9.15 of \cite{GT} (see also the monographs \cite{KS} and \cite{Mo}),
the unique solution $u$ of \hypreff{inpaint 3} and \hypreff{inpaint 4} satisfies (recall $p \geq 2$)
\[
u \in W^{1,p}_{0}(\Omega)  \cap W^{2,p}(\Omega) \psp .
\]
Next, we refer to Theorem 9.14 of \cite{GT} which yields uniform constants $c$, $C$, not depending on
such that
\begin{equation}\label{inpaint 6}
\|u\|_{W^{2,p}(\Omega)} \leq C \|g-cu\|_{L^p(\Omega)} \psp .
\end{equation}
The proposition is proved by inserting \hypreff{inpaint 2} and \hypreff{inpaint 5} 
in \hypreff{inpaint 6}. \hspace*{\fill}$\Box$\\

Recall  that $d_{\alpha\beta}$, $\alpha$, $\beta =1$, $2$, 
\[
d_{\alpha \beta}:\psp \overline{\Omega} \times \rz \times \rz^{2}\to \rz
\]
are continuous functions, in particular these functions are uniformly continuous whenever
we consider the restriction on $\overline{\Omega} \times S$ for a compact $S \subset \rz\times \rz^2$. 
Quoting the theorem of Arcela and Ascoli, it is therefore
immediate on account of \hypref{Proposition}{inpaint prop 1} that we have
condition \hypreff{intro 7} from \hypref{Assumption}{intro ass 1}.\\

It remains to justify \hypreff{intro 6} for our particular choice of $\diff$, i.e.
\[
\diff(w) := \Big(d_{\alpha\beta}\big(\cdot,I(w),\nabla I(w)\big)\Big)_{1\leq \alpha,\beta \leq 2} \psp .
\]

We first observe
\begin{proposition}\label{inpaint prop 2}
With the notation of above, the operator 
\[
I:\psp W^{1,2}(\Omega) \to W^{1,2}_0(\Omega)
\]
is continuous w.r.t.~corresponding norms of these spaces.
\end{proposition}
\emph{Proof.} We have to consider a sequence $\{w_n\}$ in $W^{1,2}(\Omega)$,
$u_n:= I(w_n)$,
\begin{equation}\label{inpaint 7}
w_n \to : w \in W^{1,2}(\Omega)\fsp\mbox{as}\msp n\to \infty\psp , \fsp u:= I(w)\psp .
\end{equation}
As above, the minimality of $u_n$ implies
\[
\iom |\nabla u_n|^2\dx \leq \lambda \int_K |w_n|^s \dx \psp .
\] 
By assumption, we have the strong convergence stated in \hypreff{inpaint 7}, thus
\[
\sup_n \int_\Omega |w_n|^s\dx < \infty
\] 
and Sobolev's inequality gives
\[
\sup_n \|u_n\|_{W^{1,2}(\Omega)} < \infty \psp .
\]
Hence, passing to a subsequence $\{\tilde{u}_n\}$, we find $\tilde{u}\in W^{1,2}_0(\Omega)$ such that
\begin{equation}\label{inpaint  8}
\tilde{u}_n \rightharpoondown \tilde{u}\fsp\mbox{in}\msp W^{1,2}(\Omega)\fsp\mbox{as}\msp
n\to \infty \psp .
\end{equation}
We now claim that we have for all $v \in W^{1,2}(\Omega)$
\begin{equation}\label{inpaint 9}
\iom |\nabla \tilde{u}|^2\dx + \lambda \int_K |\tilde{u}-w|^s \dx
\leq
\iom |\nabla v|^2\dx + \lambda \int_K |v-w|^s \dx \psp ,
\end{equation}
which immediately gives
\begin{equation}\label{inpaint 10}
\tilde{u} = I(w) = u\psp .
\end{equation}
In order to show \hypreff{inpaint 9}, we observe that this inequality holds if $\tilde{u}$ is replaced by
$\tilde{u}_n$ and that we have the lower semicontinuity
\[
\iom |\nabla \tilde{u}|^2\dx \leq \liminf_{n\to \infty} \iom |\nabla \tilde{u}_n|^2\dx \psp .
\]
In addition, \hypreff{inpaint 8} gives $\tilde{u}_n \to \tilde{u}$ in $L^q$ for any $q < \infty$,
hence our claim \hypreff{inpaint 9} and its consequence \hypreff{inpaint 10}.\\ 

Next, \hypreff{inpaint 8} implies as $n\to \infty$
\begin{equation}\label{inpaint 11}
u_n \rightharpoondown u \fsp\mbox{in}\msp W^{1,2}(\Omega)\fsp\mbox{and}\fsp
\|u_n-u\|_{L^{q}(\Omega)} \to 0 \fsp\mbox{for all}\msp q < \infty \psp .
\end{equation}

We are now going to prove
\begin{equation}\label{inpaint 12}
\nabla u_n \to \nabla u \fsp\mbox{in}\msp L^2(\Omega)\psp .
\end{equation}
As in \hypreff{inpaint 3} and \hypreff{inpaint 4} we have
\begin{eqnarray}\label{inpaint 13}
- \Delta u_n&=& g_n \fsp\mbox{on}\msp\Omega\psp ,\fsp u=0 \msp\mbox{on $\partial \Omega$}\psp , \\[2ex]
\label{inpaint 14} g_n&:=& - s |u_n-w|^{s-2} (u_n-w) \chi_K  \psp .
\end{eqnarray}
and testing the difference of 
\hypreff{inpaint 13}, \hypreff{inpaint 14} and \hypreff{inpaint 3}, \hypreff{inpaint 4} in the weak form 
with the admissible function $u_n-u$ we obtain
\[
\iom |\nabla (u_n-u)|^2\dx \leq c \iom \big(|g_n|+|g|\big)|u_n-u|\dx \psp .
\]
On account of \hypreff{inpaint 5} and \hypreff{inpaint 11} we obtain \hypreff{inpaint 12},
thus \hypref{Proposition}{inpaint prop 2}. \hspace*{\fill} $\Box$\\

Now let $w_n$, $w\in W^{1,2}(\Omega)$ s.t.~$w_n \to w$ in $W^{1,2}(\Omega)$ and define
\[
u_n := I(w_n)\psp , \fsp u:= I(w)\psp .
\]
\hypref{Proposition}{inpaint prop 2} gives as $n\to \infty$
\[
\|u_n -u\|_{W^{1,2}(\Omega)} \to 0 \psp .
\]
We then return to the proof of \hypref{Proposition}{inpaint prop 1} which gives
inequality \hypreff{inpaint 6} for the difference $u_n-u$:
\[
\|u_n-u\|_{W^{2,p}(\Omega)} \leq C \| (g_n-g)-c(u_n-u)\|_{L^p(\Omega)}\psp ,
\]
hence we have with $\alpha$ as above and as $n\to \infty$
\[
I(w_n) \to I(w) \fsp\mbox{in}\msp C^{1,\alpha}(\overline{\Omega}) \psp .
\]
We end up with
\[
\|\diff(w_n) - \diff(w)\|_{L^{\infty}(\Omega)} \to 0 \psp ,
\]
which proves \hypreff{intro 6} and thereby completes the proof of \hypref{Theorem}{diff theo 1}.\hspace*{\fill}$\Box$.\\

\bibliography{DEED_model}
\bibliographystyle{unsrt}

\vspace*{0.5cm}
\begin{tabular}{ll}
Michael Bildhauer&bibi@math.uni-sb.de\\
Martin Fuchs&fuchs@math.uni-sb.de
\end{tabular}

\vspace*{0.5cm}
\small
Department of Mathematics\\
Saarland University\\
Faculty Math.~and Computer Sci.\\
P.O.~Box 15 11 50\\
66041 Saarbr\"ucken, Germany
\end{document}